\documentclass{article}
\pdfoutput=1
\usepackage[T1]{fontenc}
\usepackage[latin9]{inputenc}
\usepackage{geometry}
\usepackage{float}
\usepackage{mathtools}
\usepackage{amsmath}
\usepackage{amssymb}
\usepackage{graphicx}
\usepackage{makecell}
\usepackage{authblk}
\makeatletter

\floatstyle{ruled}
\newfloat{algorithm}{tbp}{loa}
\providecommand{\algorithmname}{Algorithm}
\floatname{algorithm}{\protect\algorithmname}

\numberwithin{equation}{section}
\numberwithin{figure}{section}

\usepackage{color}
\usepackage{hyperref}
\usepackage{algorithm,algpseudocode}

\newcommand{\REV}[1]{{#1}}

\makeatother

\usepackage{babel}
\begin{document}
\global\long\def\ve{\varepsilon}%
\global\long\def\R{\mathbb{R}}%
\global\long\def\Rn{\mathbb{R}^{n}}%
\global\long\def\Rd{\mathbb{R}^{d}}%
\global\long\def\E{\mathbb{E}}%
\global\long\def\P{\mathbb{P}}%
\global\long\def\bx{\mathbf{x}}%
\global\long\def\vp{\varphi}%
\global\long\def\ra{\rightarrow}%
\global\long\def\smooth{C^{\infty}}%
\global\long\def\Tr{\mathrm{Tr}}%
\global\long\def\bra#1{\left\langle #1\right|}%
\global\long\def\ket#1{\left|#1\right\rangle }%
\global\long\def\Re{\mathrm{Re}}%
\global\long\def\Im{\mathrm{Im}}%
\global\long\def\bsig{\boldsymbol{\sigma}}%
\global\long\def\btau{\boldsymbol{\tau}}%
\global\long\def\bmu{\boldsymbol{\mu}}%
\global\long\def\bx{\boldsymbol{x}}%
\global\long\def\bups{\boldsymbol{\upsilon}}%
\global\long\def\bSig{\boldsymbol{\Sigma}}%
\global\long\def\bt{\boldsymbol{t}}%
\global\long\def\bs{\boldsymbol{s}}%
\global\long\def\by{\boldsymbol{y}}%
\global\long\def\brho{\boldsymbol{\rho}}%
\global\long\def\ba{\boldsymbol{a}}%
\global\long\def\bb{\boldsymbol{b}}%
\global\long\def\bz{\boldsymbol{z}}%
\global\long\def\bc{\boldsymbol{c}}%
\global\long\def\balpha{\boldsymbol{\alpha}}%
\global\long\def\bbeta{\boldsymbol{\beta}}%
\global\long\def\T{\mathrm{T}}%
\global\long\def\trip{\vert\!\vert\!\vert}%
\global\long\def\lrtrip#1{\left|\!\left|\!\left|#1\right|\!\right|\!\right|}%
\global\long\def\eps{\epsilon}%

\title{A gradient-based and determinant-free framework for \\
	fully Bayesian Gaussian process regression}
\author[1]{P. Michael Kielstra}
\author[1,2]{Michael Lindsey}
\affil[1]{University of California, Berkeley}
\affil[2]{Lawrence Berkeley National Laboratory}
\renewcommand\Affilfont{\itshape\small}
\maketitle

\begin{abstract}
Gaussian Process Regression (GPR) is widely used for inferring functions
from noisy data. GPR crucially relies on the choice of a kernel, which
might be specified in terms of a collection of hyperparameters that
must be chosen or learned. Fully Bayesian GPR seeks to infer these
kernel hyperparameters in a Bayesian sense, and the key computational
challenge in sampling from their posterior distribution is the need
for frequent determinant evaluations of large kernel matrices. This
paper introduces a gradient-based, determinant-free approach for fully
Bayesian GPR that combines a Gaussian integration trick for avoiding
the determinant with Hamiltonian Monte Carlo (HMC) sampling. Our framework
permits a matrix-free formulation and reduces the difficulty of dealing
with hyperparameter gradients to a simple automatic differentiation.
Our implementation is highly flexible and leverages GPU acceleration
with linear-scaling memory footprint. Numerical experiments demonstrate
the method's ability to scale gracefully to both high-dimensional
hyperparameter spaces and large kernel matrices.
\end{abstract}

\section{Introduction}

Gaussian process regression (GPR) \cite{Rasmussen_Williams_2006}
is a flexible framework for inferring a function on $\mathbb{R}^{d}$
from noisy observations at some large number $N$ of scattered points.
Given a positive semidefinite covariance function, or kernel, $\mathcal{K}(x,x')$,
the pointwise evaluations of a function $f$ drawn from the associated
Gaussian process are Gaussian-distributed with $\textrm{Cov}[f(x),f(x')]=\mathcal{K}(x,x')$
for any $x,x'\in\mathbb{R}^{d}$. Often one is interested in finding
the function $f$ of maximum posterior likelihood conditioned on
noisy observations, as well as quantifying the uncertainty of the
posterior.

The choice of kernel is crucial for high-quality inference. Kernels
are often grouped into families, parametrized by hyperparameters.
For example, a family of squared-exponential kernels can be parameterized
as 
\[
\mathcal{K}(x,x')=\mathcal{K}_{\alpha,\ell}(x,x')=\alpha\exp\left(-\frac{\vert x-x'\vert^{2}}{2\ell^2}\right),
\]
where $\alpha,\ell>0$ are hyperparameters that set the vertical and
horizontal scales of the Gaussian process prior. It is typical for
the user to simply fix such hyperparameters, but the goal of this
work is to infer the hyperparameters themselves in a Bayesian framework
and in particular to sample from their posterior distribution, given
our observations. This task is called fully Bayesian Gaussian process
regression \cite{pmlr-v118-lalchand20a}.

The posterior density for the hyperparameters, which from now on we
collect into a vector $\theta\in\R^{n}$, can be determined analytically
up to a normalizing constant, and therefore it is possible to draw
samples from it using Markov chain Monte Carlo (MCMC) methods. However,
the expression for the hyperparameter posterior density involves the
determinant of the $N\times N$ kernel matrix, which must be calculated
at every step in the Markov chain. This is the key computational difficulty.

We comment that direct methods \cite{Ambikasaran_Foreman-Mackey_Greengard_Hogg_ONeil_2016}
based on hierarchical matrix factorization have been applied to compute
the determinant. While this can be an excellent strategy in small
spatial dimensions ($d=1$), the scaling of such methods is in general exponential in $d$. We aim to avoid any such explicit dependence
in our approach.

Meanwhile, several previous approaches have centered around the observation
that $$\log\det A=\Tr\log A,$$ reducing the problem to one of computing
the trace of a matrix logarithm $\log A$ which can be estimated using
a randomized trace estimator \REV{\cite{10.5555/3295222.3295380,Boutsidis_Drineas_Kambadur_Kontopoulou_Zouzias_2017,Cortinovis_Kressner_2022}},
powered by linear solves or matrix-vector multiplications. Such approaches
are naturally matrix-free and therefore benefit from structure such
as sparsity of $A$. As such they have recently been used in
cutting-edge explorations of non-stationary kernel design \cite{Noack_Krishnan_Risser_Reyes_2023,NoackLuoRisser2024}
where sparsity is critical for scalability\@.

Nonetheless, \REV{these approaches have several drawbacks.} First, the log-determinant is not computed exactly
in this framework: \REV{the matrix logarithm itself is usully approximated rather than exactly evaluated, and the trace is calculated using a randomized estimator. The resulting approximation suffers from} a slowly converging
Monte Carlo error rate and introduces a bias
into the posterior sampler that is difficult to control \emph{a priori}.
Relatedly, the cost of driving this error to negligibility is severe,
since several linear solves or matrix-vector multiplications may be
required for each \REV{estimation of the trace. In turn, the entire procedure
is wrapped in an expensive Markov chain sampler.}

In \cite{determinantfree} an approach was proposed which avoids the
determinant entirely, sidestepping the quandary of determining the
numerical tolerance for the log-determinant approximation. The approach
is based on the observation that for a positive definite matrix $A$,
\[
\left[\det A\right]^{-1/2}\propto\int_{\mathbb{R}^{N}}e^{-\frac{1}{2}\varphi^{\top}A\varphi}\,d\varphi.
\]
This observation allows us to introduce an auxiliary variable $\varphi$
which is Gaussian conditioned on $\theta$ and which, when integrated
out, yields the desired target distribution for $\theta$. We can
then sample from the joint distribution for $(\theta,\vp)$ and simply
discard our samples of $\varphi$, leaving behind samples of $\theta$
which do not require computing determinants at all.

This approach can also be viewed as being rooted in the longer history
of lattice \REV{quantum chromodynamics} (LQCD) calculations \cite{Fucito_Marinari_Parisi_Rebbi_1981,Clark_Kennedy_2007,Gattringer_Lang_2010},
which require samples from a probability density with formally similar
structure. In that setting, the introduction of the auxiliary variable
is widely known as the `pseudofermion trick.' To carry out the pseudofermion
trick in the setting of fully Bayesian GPR, some additional techniques
for the \REV{sampling} of $\vp$ (conditioned on $\theta$) are required,
beyond what is usually required in LQCD, as we shall elaborate below.
In particular, rational approximation techniques \cite{doi:10.1137/070700607, HOFREITHER2020332}
must be used to perform matrix-vector multiplications by the kernel
matrix square root. In addition to the setting of fully Bayesian GPR
\cite{determinantfree}, such rational approximation techniques have
been used to optimize the kernel hyperparameters in an MLE or MAP
sense \cite{10.5555/3495724.3497591}.

In spite of its advantages, the work \cite{determinantfree} on determinant-free
fully Bayesian GPR does not leverage gradients in its MCMC routine.
When the number $n$ of hyperparameters becomes large---as advanced
approaches to flexible kernel design increasingly demand \cite{Noack_Krishnan_Risser_Reyes_2023,NoackLuoRisser2024}---the
advantage of using gradients in the MCMC routine becomes more dramatic.
We comment that a recent work \cite{latz2024deepgaussianprocesspriors},
extending the aforementioned work \cite{determinantfree}, has implemented
a preconditioned Crank-Nicolson sampler \cite{Cotter_Roberts_Stuart_White_2013,Hairer_Stuart_Vollmer_2014},
which is motivated by the fact that their hyperparameters arise from
an increasingly fine spatial discretization. Although such an approach
can work well for such high-dimensional discretizations, it can struggle
when the target distribution is far from Gaussian, or when the covariance
structure of the target distribution cannot be compressed tractably.
The flexibility and widespread success of gradient-based approaches
across a variety of sampling tasks motivates our work.

We pursue an approach inspired by the standard practice of the LQCD
community, which couples the pseudofermion trick with Hybrid/Hamiltonian
Monte Carlo (HMC) \cite{Duane_Kennedy_Pendleton_Roweth_1987} sampling
for the physical variable $\theta$. (In fact, HMC was invented for
this purpose.) The strategy for computing the ``force \REV{term}'' required by
HMC (cf. Section \ref{subsec:Forces}) relies on only a few computational
primitives involving the kernel matrix $A(\theta)$, which can all
be viewed as matrix-free in the sense that the kernel matrix is only used for matrix-vector multiplications and never explicitly formed. In particular, we make use of a formulation
in which the computation of gradients can be reduced to the primitive
operation of differentiation through quadratic forms involving the
kernel matrix, i.e., the computation of $\nabla_{\theta}\left[z^{\top}A(\theta)z\right]$
for arbitrary $z\in\R^{N}$.

Our computational primitives can all be implemented conveniently using
the KeOps library \cite{JMLR:v22:20-275}, which runs efficiently
on GPUs and supports flexible kernel design as well as efficient automatic
differentiation for quantities such as $z^{\top}A(\theta)z$. Moreover,
without relying on any sparsity or low-rank structure, the memory footprint
of all required operations is only $O(N)$. While the cost scaling
is generically quadratic in $N$ for dense unstructured kernel matrices,
in practice the hardware acceleration offered by KeOps is so dramatic
that our approach may exceed the performance of alternative
implementations that are more aware of the structure of $A$. Meanwhile, the completely generic and flexible quality
of the implementation offers the potential for widespread use. Finally,
as we shall see, the additional cost of computing gradients is essentially
negligible for large kernel matrices, as it is dominated by the cost
of the linear solves required to implement the forces. \REV{The theoretical} advantages \REV{of gradient-based methods}
over gradient-free Random Walk Metropolis \REV{are well-understood \cite{Chewi_2024}, and our experiments also demonstrate advantages in practice}. These advantages
notwithstanding, we \REV{note} that our implementation choices are not
necessarily essential. As matrix-free approaches, our algorithms could
benefit from sparsity or efficient routines for matrix-vector multiplications
\cite{greengard_equispaced_2023,barnett_uniform_2023,kielstra2024gaussianprocessregressionloglinear}.

To our surprise, the work \cite{determinantfree} on determinant-free
fully Bayesian GPR seems to have received little attention in the
literature, besides the recent aforementioned extension \cite{latz2024deepgaussianprocesspriors}.
Various other works \cite{NEURIPS2018_4172f310,pmlr-v130-rossi21a}
have avoided the difficulty of dealing with the determinant by adopting
a stochastic gradient HMC framework \cite{pmlr-v32-cheni14}, in which
the collection of observations are broken randomly into mini-batches.
This procedure reduces the size of the kernel matrix operations at
the expense of introducing some bias. (We remark that the mini-batches
could themselves be scaled to be larger by applying the pseudofermion
trick to avoid computing their determinants.) Meanwhile, HMC has been
applied more straightforwardly (without either randomized trace estimation
or the pseudofermion trick), e.g., in \cite{pmlr-v51-heinonen16},
with exact Metropolis correction, albeit with scaling at least $O(N^3)$ per step due to direct computation of determinants.
To our knowledge, it seems not to be appreciated that many of the
advantages of these preceding works can be combined, achieving asymptotically
unbiased HMC with $O(N^2)$ operations per step.

\subsection{Acknowledgements}

This work used Jetstream2 at Indiana University through allocation
MTH240042 from the Advanced Cyberinfrastructure Coordination Ecosystem:
Services \& Support (ACCESS) program, which is supported by National
Science Foundation grants \#2138259, \#2138286, \#2138307, \#2137603,
and \#2138296. For more on ACCESS and Jetstream2 see \cite{Boerner_Deems_Furlani_Knuth_Towns_2023}
and \cite{Hancock_Fischer_Lowe}. M.L. was partially supported by
a Sloan Research Fellowship as well as the Applied Mathematics Program
of the US Department of Energy (DOE) Office of Advanced Scientific
Computing Research under contract number DE-AC02-05CH11231.

\section{Preliminaries \label{sec:prelim}}

We are given noisy measurements $y_{i}\in\R$, $i=1,\ldots,N$, of
an unknown function $f:\R^{d}\ra\R$, on a set of scattered points
$x_{i}\in\R^{d}$. More concretely, we are given values
\begin{equation}
y_{i}=f(x_{i})+\eps_{i},\label{eq:model}
\end{equation}
 where $\eps_{i}$ are independently distributed Gaussian noise terms.

We consider a Gaussian process prior over the function $f$, \REV{one} induced
by a positive semidefinite kernel $\mathcal{K}_{\theta}(x,x')$, where
the kernel itself is parameterized by hyperparameters $\theta\in\R^{n}$
that we aim to infer. We can also accommodate uncertainty in the distribution
of the noise terms $\eps_{i}$ by allowing their independent distributions
to depend on the unknown parameters $\theta$. To wit, letting $\sigma_{\theta}^{2}(x)>0$
define a spatially varying noise variance, we assume 
\[
\eps_{i}\sim\mathcal{N}(0,\sigma_{\theta}^{2}(x_{i})).
\]
Finally, let $p(\theta)$ denote a prior distribution over the uncertain
hyperparameters. The prior distribution for $f$ depends on $\theta$: while the overall density $p(f)$ might not be Gaussian, $p(f\vert\theta)$ is.

Under these conditions \cite{Rasmussen_Williams_2006}, the posterior
distribution $P(\theta)=p(\theta\vert y)$\footnote{Since we will view $y$ as fixed in our discussion, for simplicity
we opt to omit it from the notation for our target distribution $P(\theta)$.} over \REV{the vector $\theta$ of} hyperparameters, given our observations $y$, is given by
\begin{equation}
P(\theta)\propto\left|A(\theta)\right|^{-1/2}e^{-\frac{1}{2}y^{\top}A(\theta)^{-1}y}\,p(\theta),\label{eq:P}
\end{equation}
 where $A(\theta)\in\R^{N\times N}$ is a positive definite matrix
with entries defined by 
\begin{equation}
A_{ij}(\theta)=\sigma_{\theta}^{2}(x_{i})\:\delta_{ij}+\mathcal{K}_{\theta}(x_{i},x_{j}),\label{eq:kerdef}
\end{equation}
 and $\vert\,\cdot\,\vert$ indicates the matrix determinant. \REV{We also define \begin{equation}\label{eq:kernmat}
K_{ij}(\theta)=\mathcal{K}_{\theta}(x_{i},x_{j}),
\end{equation} i.e., the kernel matrix, which is useful in our discussion of preconditioning.}  Usually,
we take $\sigma_{\theta}^{2}(x)\equiv\sigma^{2}$ to be constant,
and in this case we can think of $\sigma$ as either fixed or as part
of the collection of hyperparameters $\theta$ to be inferred.

The task of \emph{fully Bayesian Gaussian process regression }\cite{pmlr-v118-lalchand20a}
is to draw samples defined by this posterior density $P(\theta)$.
The key computational difficulty is imposed by the determinant, naive
evaluation of which has a cost that scales as $O(N^{3})$. 

\section{Methodology}

Inspired by the so-called \emph{pseudofermion trick} that is widely
used in lattice QCD \cite{Fucito_Marinari_Parisi_Rebbi_1981,Clark_Kennedy_2007,Gattringer_Lang_2010}
and, more recently, in condensed matter physics \cite{Scalettar_Scalapino_Sugar_Toussaint_1987,Lunts_Albergo_Lindsey_2023},
we propose to deal with the determinant by viewing it as a Gaussian
integral: 
\[
\left[\det A(\theta)\right]^{-1/2}\propto\int_{\mathbb{R}^{N}}e^{-\frac{1}{2}\varphi^{\top}A(\theta)\varphi}\,d\varphi.
\]
 Accordingly, we can view $P(\theta)$ as the marginal of the extended
density $P(\theta,\vp)$ involving an auxiliary variable $\vp\in\R^{N}$:
\begin{equation}
\label{eq:extendeddensity}
P(\theta,\varphi)\propto e^{-\left[S(\theta)+\frac{1}{2}y^{\top}A(\theta)^{-1}y+\frac{1}{2}\varphi^{\top}A(\theta)\varphi\right]},
\end{equation}
 where we have defined $S(\theta)$ such that 
\[
p(\theta)\propto e^{-S(\theta)}.
\]
 Then provided that we can draw samples $(\theta,\vp)$ from $P(\theta,\vp)$,
we obtain samples from $P(\theta)$ simply by omitting the auxiliary
variable.

As in the physics literature, we will propose to sample from $P(\theta,\vp)$
using two key routines:
\begin{itemize}
\item [(1)] Exact Gibbs \REV{updating} of the $\vp$, given fixed $\theta$.
\item [(2)] Random Walk Metropolis or Hamiltonian Monte Carlo sampling
of $\theta$, given fixed $\vp$.
\end{itemize}
To achieve step (1), we will see that we encounter a significant,
but \REV{not insurmountable}, complication \REV{that does not appear in} the Gibbs \REV{update} of the
pseudofermion field in the physical literature. Step (2) may be achieved
through many possible methods. We consider standard Random Walk Metropolis \cite{Betancourt_2018} or
one of two Hamiltonian Monte Carlo approaches: a
standard technique and a new approach based on implicit integrators.

\subsection{\REV{Updating} the auxiliary variable}\label{subsec:varphi}

The conditional density for $\vp$ given $\theta$ is given by
\[
P(\vp\:\vert\,\theta)\propto e^{-\frac{1}{2}\vp^{\top}A(\theta)\vp},
\]
 i.e., $\vp\:\vert\,\theta$ is distributed as $\mathcal{N}(0,A^{-1}(\theta))$.
In order to exactly sample $\vp\:\vert\,\theta$, we can sample a
standard normal random variable $\xi\sim\mathcal{N}(0,I_{N})$ and
then form $\vp=A^{-1/2}(\theta)\,\xi$, furnishing\footnote{By contrast, in the physical literature, the `pseudofermion field'
$\vp$ is typically sampled from a Gaussian distribution with covariance
$D(\theta)D(\theta)^{*}$, where $D(\theta)$ indicates the `fermion
matrix' (and $\theta$ plays the role of the bosonic field). In these
settings $D(\theta)$ is typically a sparse matrix and the pseudofermion
field can be sampled straightforwardly by the sparse matrix multiplication
$D(\theta)\,\xi$.} an exact sample from $P(\vp\,\vert\,\theta)$.

We can approximate the matrix-vector multiplication $A^{-1/2}(\theta)\,\xi$
by a collection of linear solves by shifts of $A(\theta)$, using
the the pole expansion \cite{doi:10.1137/070700607} 
\begin{equation}
A^{-1/2}(\theta)\approx\sum_{p=1}^{N_{p}}w_{p}\,(A(\theta)+\lambda_{p}I_{N})^{-1},\label{eq:pex}
\end{equation}
 in which the weights $w_{p}$ are real and the shifts $\lambda_{p}\geq0$
are nonnegative. Moreover, the approximation is rapidly converging
in the number $N_{p}$ of poles, and the choice of $N_{p}$ required
for fixed accuracy scales only logarithmically with the condition
number of $A(\theta)$. See Appendix \ref{app:hht} for further concrete
details. Note that since the shifts are nonnegative, each of the linear
solves demanded by the pole expansion is positive definite.

\REV{Certain preconditioners for $A(\theta)$ can be adapted to the shifted solves as well.  Specifically, in some experiments we will consider Nystr\"om preconditioners, which are constructed from an approximate low-rank factorization of the kernel matrix $K(\theta)$~\eqref{eq:kernmat}. In particular we use the randomized Nystr\"{o}m approach of \cite{frangella_tropp_udell}, though it is also possible to use approaches based on column selection (see, e.g., \cite{chen2023randomly,fornace2024columnrowsubsetselection} for reviews and recent examples).
Since $A(\theta) - K(\theta)$ is a diagonal matrix, this low-rank factorization can be used to find an approximation to $A(\theta)^{-1}$ using the Woodbury matrix identity \cite{frangella_tropp_udell}.  The shift terms $\lambda_p I_N$ appearing in~\eqref{eq:pex} only affect the diagonal, so the Woodbury identity can also be applied to form an effective preconditioner for each of the solves necessitated by~\eqref{eq:pex}, recycling the same low-rank factorization of $K(\theta)$ for each term.}

\subsection{Sampling the original variable with HMC}\label{subsec:hmc}

Exact sampling of $\theta\,\vert\,\vp$ according to the conditional density
\begin{equation}
P(\theta\,\vert\,\vp)\propto e^{-U_{\vp}(\theta)},\quad\text{where }\ U_{\vp}(\theta):=S(\theta)+\frac{1}{2}y^{\top}A(\theta)^{-1}y+\frac{1}{2}\varphi^{\top}A(\theta)\varphi,\label{eq:U}
\end{equation}
 is not computationally feasible. We will approach this task with
a Markov chain Monte Carlo (MCMC) update, specifically Hamiltonian
(alternatively, hybrid) Monte Carlo (HMC) \cite{Betancourt_2018}.
Although other generic MCMC frameworks (such as overdamped Langevin
dynamics \cite{Girolami_Calderhead_2011}) could be considered, the
HMC framework in particular allows for convenient Metropolis correction
of proposals based on implicit integrators for the Hamiltonian dynamics,
which we shall explore in more detail below.

HMC is based on the introduction of an auxiliary momentum degree of
freedom $\pi\in\R^{n}$ which is conjugate to $\theta$. We aim to
sample $\theta,\pi\sim P_{\vp}(\theta,\pi)$ defined by 
\[
P_{\vp}(\theta,\pi)\propto P(\theta\,\vert\,\vp)\,e^{-\frac{1}{2}\pi^{\top}M^{-1}\pi},
\]
 where $M\succ0$ is a positive definite `mass matrix.' \REV{We set $M=I$ in our own experiments, but, following \cite{Betancourt_2018}, other implementations might consider taking $M^{-1}$ to be the covariance matrix of $\theta$, estimated in practice using a short preliminary run with $M=I$. Choosing an effective mass matrix is essentially equivalent to adopting a linear change of coordinates that transforms the covariance matrix to the identity.}  Then samples
$\theta\sim P(\theta\,\vert\,\vp)$ can be obtained simply by omission
of the conjugate variable. We can \REV{write} 
\[
P_{\vp}(\theta,\pi)\propto e^{-H_{\vp}(\theta,\pi)}
\]
 where 
\[
H_{\vp}(\theta,\pi):=U_{\vp}(\theta)+\frac{1}{2}\pi^{\top}M^{-1}\pi
\]
 can be viewed as a Hamiltonian.

Then HMC is based on the induced Hamiltonian dynamics 

\begin{equation}
\begin{cases}
\dot{\theta}(t)=M^{-1}\pi(t)\\
\dot{\pi}(t)=-F_{\vp}(\theta(t)),
\end{cases}\label{eq:Hamdyn}
\end{equation}
 in which the `force' $F_{\vp}(\theta)$ is defined by 
\[
F_{\vp}(\theta):=\nabla_{\theta}U_{\vp}(\theta).
\]
 We will let $\mathcal{S}_{\vp,\Delta t}$ denote the discrete-time
flow map of a \emph{symplectic} integrator for these dynamics, for
a time step of $\Delta t$. Below we will discuss several options
for the symplectic integrator. In addition to the step size $\Delta t$,
we must also choose the number $N_{\mathrm{int}}$ of integration
steps performed within each update of our Markov chain.

To generate a Markov chain with the target invariant distribution
$P_{\vp}(\theta,\pi)$, HMC involves the repetition of the following
steps:
\begin{itemize}
\item [(1)] Sample a new momentum variable $\pi\sim\mathcal{N}(0,M)$.
\item [(2)] Evolve according to discretized Hamiltonian dynamics: 
\[
(\theta',\pi')=\left[\mathcal{S}_{\vp,\Delta t}\right]^{N_{\mathrm{int}}}(\theta,\pi).
\]
\item [(3)] Accept the update, i.e., replace $(\theta,\pi)\leftarrow(\theta',\pi')$,
with probability 
\[
\min\left(1,e^{H_{\vp}(\theta,\pi)-H_{\vp}(\theta',\pi')}\right).
\]
 Otherwise do not change $(\theta,\pi)$.
\end{itemize}

\REV{We usually take $N_{\mathrm{int}}>1$, and note that only the final result after all $N_{\mathrm{int}}$ steps of numerical integration is considered as a new candidate sample.}  In the following subsections we will discuss the practical computation
of the force $F_{\vp}$, which imposes the main bottleneck
on the algorithm, as well as different choices of symplectic integrator.

\subsubsection{Forces\label{subsec:Forces}}

Direct computation gives the expression for the force: 
\[
F_{\vp}(\theta)=\nabla_{\theta}S(\theta)-\frac{1}{2}y^{\top}A(\theta)^{-1}\frac{\partial A(\theta)}{\partial\theta}A(\theta)^{-1}y+\frac{1}{2}\vp^{\top}\frac{\partial A(\theta)}{\partial\theta}\vp.
\]
 However, it is inefficient to form the tensor $\frac{\partial A(\theta)}{\partial\theta}$
of partial derivatives of $A(\theta)$ with respect to the hyperparameters
$\theta$ within each force evaluation.

An alternative approach might be to automatically differentiate an evaluation of
 $U_{\vp}(\theta)$ following (\ref{eq:U}). However,
note that such an evaluation requires solving the linear system $A(\theta)x=y$.
Naively differentiating through an iterative solver for this positive
definite linear system is also a computational cost that we prefer
to avoid.

Importantly, our capacity to evaluate the force hinges only on our
capacity to perform the following subroutines:
\begin{itemize}
\item [(1)] Perform linear solves of the form $A(\theta)\,x=y$, where
$\theta$ is fixed and $x$ is unknown.
\item [(2)] Evaluate $\nabla_{\theta}S$.
\item [(3)] Evaluate the gradient with respect to $\theta$ of a quadratic
form in $A(\theta)$, i.e., evaluate $\nabla_{\theta}\left[z^{\top}A(\theta)z\right]$
for arbitrary $z\in\R^{N}$ independent of $\theta$.
\end{itemize}
Indeed, we can rewrite $F_{\vp}(\theta)$ as 
\[
F_{\vp}(\theta)=f_{\vp}(\theta,x_{\theta}),
\]
 where we define $x_{\theta}$ as the solution of the positive definite
linear system 
\begin{equation}
A(\theta)\,x_{\theta}=y\label{eq:xtheta}
\end{equation}
 and $f_{\vp}(\theta,x)$ as the gradient 
\[
f_{\vp}(\theta,x)=\nabla_{\theta}u_{\vp}(\theta,x)
\]
 of the auxiliary function 
\begin{equation}
u_{\vp}(\theta,x):=S(\theta)-\frac{1}{2}x^{\top}A(\theta)x+\frac{1}{2}\vp^{\top}A(\theta)\vp.\label{eq:u}
\end{equation}

Then to evaluate $F_{\vp}(\theta)$, we can first form $x=x_{\theta}$,
relying only on capacity (1), and then --- with the dependence of
$x=x_{\theta}$ on $\theta$ frozen --- differentiate through the
evaluation (\ref{eq:u}) of $u_{\vp}(\theta,x)$, relying only on
capacities (2) and (3). \REV{Note that the linear solves~\eqref{eq:xtheta} can be preconditioned, e.g., following the discussion at the end of Section~\ref{subsec:varphi}.}

\subsubsection{Explicit integrator\label{subsec:explicit}}

The standard choice of symplectic integrator for HMC is the St\"ormer-Verlet
or `leapfrog' integrator \cite{Betancourt_2018}, which is second-order
accurate. Given $(\theta_{0},\pi_{0})$, the map $\mathcal{S}_{\vp,\Delta t}$
is defined by $(\theta_{0},\pi_{0})\mapsto(\theta_{1},\pi_{1})$,
where the image is determined by the algorithmic steps: 
\begin{align*}
 & \theta_{1/2}\leftarrow\theta_{0}+\frac{\Delta t}{2}\,M^{-1}\pi_{0}\\
 & \pi_{1}\leftarrow\pi_{1}-(\Delta t)\,F_{\vp}(\theta_{1/2})\\
 & \theta_{1}\leftarrow\theta_{1/2}+\frac{\Delta t}{2}\,M^{-1}\pi_{1}.
\end{align*}
 Note that in fact the repeated application $[\mathcal{S}_{\vp,\Delta t}]^{N_{\mathrm{int}}}$
requires only $N_{\mathrm{int}}+1$ inversions of the mass matrix
and $N_{\mathrm{int}}$ force evaluations, since the last and first
stages of the update can be combined into a single step for all but
the first and last of the $N_{\mathrm{int}}$ time steps.

\subsubsection{Implicit integrator\label{subsec:implicit}}

Observe that even in the implementation of an explicit integrator
such as the leapfrog integrator, we must solve the linear system (\ref{eq:xtheta})
once per integration step. The linear solver itself, for large systems,
will be iterative.

Meanwhile, the typical case for using an explicit integrator is that
an iterative approach for solving the integration scheme is not required.
Therefore we are motivated to consider implicit symplectic integrators.
If we can devise an iterative scheme that simultaneously approaches a solution to
the nonlinear equations defining the implicit integrator \emph{as well
as} the linear system defining the force, then the typical case against
using an implicit integrator may dissipate. Meanwhile, if the implicit
integrator allows for a larger timestep $\Delta t$, then we may \REV{gain an overall computational advantage by reducing the mixing time of the Markov chain sampler.}

Specifically, we will consider the simplest possibility, the 1-stage
Gauss-Legendre or `implicit midpoint' scheme. All Gauss-Legendre integrators
are symplectic \cite{Sanz-Serna_1988}.

For our Hamiltonian dynamics (\ref{eq:Hamdyn}), the implicit midpoint
integrator is defined by the following equations in the unknowns $(\theta_1,\pi_1,x)$:
\begin{equation}
\begin{cases}
\theta_1=\theta_{0}+(\Delta t)\,M^{-1}\left[\frac{\pi_{0}+\pi_1}{2}\right],\\
\pi_1=\pi_{0}-(\Delta t)\,f_{\vp}\left(\frac{\theta_{0}+\theta_1}{2},x\right),\\
A\left(\frac{\theta_{0}+\theta_1}{2}\right)\,x=y.
\end{cases}\label{eq:implicitsystem}
\end{equation}
 (We do still have to solve $A(\theta_1)x=y$, which is a different
equation, for the Metropolis-Hastings routine, but we only need to
solve it a single time once we have \REV{found $\theta_1$}.) Assuming
the existence of a unique solution $(\theta_1,\pi_1,x)$ to these equations,
the integrator map $\mathcal{S}_{\vp,\Delta t}$ is defined by the
update $(\theta_{0},\pi_{0})\mapsto(\theta_1,\pi_1)$.

To solve these equations, we introduce the map $\Phi:\R^{n}\times\R^{n}\times\R^{N}\ra\R^{n}\times\R^{n}\times\R^{N}$
defined by 
\begin{equation}
\Phi\left(\begin{array}{c}
\theta\\
\pi\\
x
\end{array}\right)=\left(\begin{array}{c}
\theta_0+(\Delta t)\,M^{-1}\left[\frac{\pi_{0}+\pi}{2}\right]\\
\pi_0-(\Delta t)\,f_{\vp}\left(\frac{\theta_{0}+\theta}{2},x\right)\\
x+R(\theta_{0})\,\left[y-A\left(\frac{\theta_{0}+\theta}{2}\right)x\right]
\end{array}\right),\label{eq:Phi}
\end{equation}
where \REV{$R(\theta_{0})$} is a preconditioner, possibly chosen dependent
on $\theta_{0}$, for the linear system $A(\theta_{0})\,x=y$.

We seek a fixed point of $\Phi$ (\ref{eq:Phi}), which is equivalent to a solution to (\ref{eq:implicitsystem}). Fixed-point iteration
by $\Phi$ can be viewed as a fusion of standard Picard iteration
for the implicit midpoint method with preconditioned gradient descent
for the positive definite linear system appearing in (\ref{eq:implicitsystem}).

In practice, we will accelerate the fixed-point iteration by $\Phi$
using \emph{Anderson acceleration} \cite{Anderson_1965, CHEN2024116077}. In fact, ignoring the nonlinear
dependence on $\theta$, the behavior of Anderson acceleration is
well-understood \cite{Walker_Ni_2011} when applied to the linear
fixed-point iteration 
\[
x\mapsto x+R[y-Ax]
\]
 for solving the linear system $Ax=y$. Indeed the performance of
Anderson acceleration in this simple context is quite closely connected
to that of GMRES. Thus our approach can be viewed as a robust extension
of GMRES for the linear system, which simultaneously converges to a solution of the
nonlinear equations required by the implicit integrator. More general
theory connecting the convergence of Anderson acceleration to that
of quasi-Newton methods has been explored recently in \cite{feng2024convergenceanalysisalternatingandersonpicard}.

The use of a preconditioner $R$ is supported by \cite[Corollary 3.5]{Walker_Ni_2011}
in particular, which connects the solution of this exact fixed-point
equation to the solution of $RAx=Ry$ using GMRES.  \REV{A requirement for the non-accelerated fixed-point iteration to converge, and therefore a requirement for the Anderson-accelerated iteration to reliably converge as well, is that the spectral radius $\rho(RA)$ is bounded by 1.  We therefore have two goals in choosing our preconditioner.  First and foremost, it should simply rescale the matrix so that $RA$ is sufficiently small. Ideally it should also approximate $A^{-1}$ sufficiently well to accelerate a hypothetical GMRES solve of $RAx=Ry$.}

\REV{The first goal can be fulfilled by a simple rescaling: we estimate $c := \rho( A )$ with a few iterations of the power method and then set $R=\frac{1}{c}I$. The second goal can be fulfilled by a more sophisticated construction of an approximate inverse. Since $A$ is a diagonal update to the kernel matrix $K$ (cf. \eqref{eq:kerdef}-\eqref{eq:kernmat}) we can construct a low-rank approximation of $K$ and in turn an approximation of $A^{-1}$ via the Woodbury identity \cite{frangella_tropp_udell}. When pursuing such an approach, we use the randomized Nystr\"{o}m approach of \cite{frangella_tropp_udell}. We can combine this approach with a rescaling step, first calculating an initial preconditioner $\tilde{R}$ using the Nystr\"om method and then, using a power method to find $\tilde{c} \coloneqq \rho( \tilde{R}A) $, setting $R=\frac{1}{\tilde{c}}\tilde{R}$.  In practice, in our experiments, we found that this did not meaningfully improve or accelerate convergence.}

\REV{Determining sufficient conditions for the full non-linear iteration to converge is more complicated and falls outside the scope of this work. More broadly, the convergence of Anderson acceleration and its variants is a well-studied and deep problem \cite{Toth_Kelley_2015, Evans_Pollock_Rebholz_Xiao_2020, Feng_Laiu_Strohmer_2024}.  The existence of a locally unique solution in the small $\Delta t$ regime is guaranteed by the implicit function theorem, but it is difficult to ascertain the extent of this regime. Meanwhile, we comment that other algorithms can also be used to solve \eqref{eq:implicitsystem}.  We experimented with Jacobian-free Newton-Krylov methods \cite{KNOLL2004357}, as well as a ``CG-within-Anderson'' approach in which the linear system was solved for $x$ to high accuracy within an Anderson-accelerated fixed-point iteration for $(\theta,\pi)$. In practice, we found that pure Anderson acceleration was the most effective scheme.}

Further background on the practicalities of Anderson acceleration,
as well as its use in our setting, are provided in Appendix
\ref{app:aa}.

\subsection{Summary}

\begin{algorithm}
\caption{Pseudocode for full Bayesian posterior hyperparameter sampler}

\begin{algorithmic}[1]
\Require{Step size $\Delta t >0$; integration count $N_{\mathrm{int}}$; symplectic integrator $\mathcal{S}_{\vp, \Delta t}$; positive definite $N\times N$ mass matrix $M$; initialization $\theta^{(0)}$; preconditioner $P(\theta)$ }
\For{$t = 0,\ldots, T-1$}
\State{Sample $\xi \sim \mathcal{N}(0,I_N)$}
\State{Form $\vp \approx A(\theta)^{-1/2} \xi$  using preconditioned linear solves via the pole expansion \eqref{eq:pex} }
\State{Sample $\pi \sim \mathcal{N}(0,M)$}
\State{Form $(\theta',\pi')=\left[\mathcal{S}_{\vp, \Delta t}\right]^{N_{\mathrm{int}}}(\theta^{(t)},\pi)$}
\State{Set $\alpha \leftarrow \min\left(1,e^{H_{\vp}(\theta^{(t)},\pi)-H_{\vp}(\theta',\pi')}\right)$, and draw $X$ uniformly randomly from $[0,1]$}
\If{$X \leq \alpha$}
\State{Set $\theta^{(t+1)} \leftarrow \theta'$}
\Else
\State{Set $\theta^{(t+1)} \leftarrow \theta^{(t)}$}
\EndIf
\EndFor
\State{\Return Markov chain trajectory $\theta^{(0)}, \ldots, \theta^{(T)}$}
\end{algorithmic}\label{alg:sampler}
\end{algorithm}

In Algorithm \ref{alg:sampler} we synthesize the steps of our sampler
for $P(\theta)$ into an abstract pseudocode. We refer the reader
to the main body for the implementation details of the pole expansion
(\ref{eq:pex}), as well as of the symplectic integrator $\mathcal{S}_{\vp,\Delta t}$,
for which we introduced two possible choices in Sections \ref{subsec:explicit}
and \ref{subsec:implicit} above.

Note that in practice, the sampler may be run for a burn-in period,
from which all samples are discarded in downstream analysis. In practice,
the HMC hyperparameters (i.e., the step size $\Delta t$, the number
$N_{\mathrm{int}}$ of integration timesteps per accept-reject step,
and the mass matrix $M$) may all be determined on-the-fly during
the burn-in period, as is standard practice in many statistical applications,
cf. \cite{Stan}.

\section{Numerical experiments}

First we present a flexible family of kernels with many hyperparameters.
We will assume that our scattered points $x_{1},\dots,x_N$ lie
in the reference box $[-1,1]^{d}$. In general, this can be achieved
by suitably shifting and scaling the domain.

Let $T_{n}(x)$ be the $n$-th Chebyshev polynomial of the first kind.
Recall that $d$ denotes the spatial dimension (i.e., the scattered
points are elements of $\R^{d}$). We identify \REV{a block of} \REV{entries in} the vector $\theta$ of hyperparameters
with a $d$-tensor $\Theta\in\mathbb{R}^{N_{\textrm{cheb}}\times N_{\textrm{cheb}}\times\cdots\times N_{\textrm{cheb}}}$,
where $N_{\mathrm{cheb}}$ denotes the order of our Chebyshev expansion.
In general we will use superscripts for component indices, writing
$x=(x^{1},\ldots,x^{d})\in\R^{d}$ for a general point and $x_{i}=(x_{i}^{1},\ldots,x_{i}^{d})$
for our scattered points.

Then we define
\[
C_{\theta}(x)=\sum_{i_{1},\dots,i_{d}=0}^{N_{\textrm{cheb}}-1}\Theta_{i_{1}\cdots i_{d}}T_{i_{1}}(x^{1})\cdots T_{i_{d}}(x^{d}).
\]

We consider the squared-exponential kernel with non-stationary vertical
scale: 
\[
\mathcal{K}_{\theta}(x,y)=\exp(C_{\theta}(x))\exp(C_{\theta}(y))\exp\left(-\frac{\vert x-y\vert^{2}}{2\ell^2}\right).
\]
\REV{This can easily be verified to be positive definite and therefore a valid kernel function.} We will take the noise variance $\sigma_{\theta}^{2}(x)\equiv\sigma^{2}$
to be constant (cf. (\ref{eq:kerdef})) and sometimes view $\sigma$
\REV{as dependent on} an additional element of $\theta$ to be inferred. (By default
we view $\sigma$ as fixed and will comment below when we attempt
to infer it.) \REV{We will treat the horizontal length scale $2\ell^2$ the same way, and observe that we} can
create as many \REV{additional} hyperparameters as we wish by taking $N_{\textrm{cheb}}$
larger.  \REV{Though we treat $\ell$ and $\sigma$ as position-independent, we comment that more generally it is possible to infer spatial dependencies in these quantities.}

We adopt a flat prior $p(\theta)\equiv1$ for the hyperparameters
and take the mass matrix in HMC to be the identity. \REV{When, as indicated above, we aim to infer other hyperparameters determining $\ell$ and $\sigma$, we adopt a flat prior for these hyperparameters as well.
Whenever we need to explicitly solve a linear system (e.g., when we compute $x=A(\theta)^{-1}y$ within in our explicit leapfrog integrator), we converge the solver up to a relative residual error tolerance of $10^{-6}$. To compute inverse matrix square roots, we use $N_p=15$ terms in the pole expansion, as discussed in Appendix \ref{app:hht}. We} always initialize $\theta$
to have all elements equal to $0.01$.

We carry out three sets of numerical experiments. First, we verify
the correctness of our `pseudofermion' methods in a simple setting with
few hyperparameters where $P(\theta)$ can be \REV{determined} by
quadrature. \REV{Throughout we use the term `pseudofermion' to indicate that we are sampling the extended density $P(\theta ,\vp)$ \eqref{eq:extendeddensity} in which $\vp$ plays the role of the `pseudofermion field' in the physical nomenclature. We use this language in order to distinguish from approaches that sample the marginal $P(\theta)$ \eqref{eq:P} directly using exact determinant computation, since we will show some comparisons between our pseudofermion-based samplers and analogous determinant-based samplers.}

\REV{Indeed, in our second set of experiments,} we compare the performance of three proposal 
mechanisms (Random Walk Metropolis and \REV{Hamiltonian Monte Carlo with both explicit and implicit integrators), both for the determinant-based and pseudofermion-based targets.  Our HMC implementation is described in Section \ref{subsec:hmc}. Our RWM implementation largely follows \cite{Betancourt_2018}: given an initial vector of hyperparameter values $\theta$, we generate a new proposal $\theta' \sim \mathcal{N}(\theta, (\Delta t)^2 I)$.  We accept this with probability $\min \left(1, e^{U_\varphi(\theta)-U_\varphi(\theta')}\right)$, otherwise leaving $\theta$ unchanged.  The sampling of $\varphi$ proceeds exactly as described in Section~\ref{subsec:varphi}. The RWM sampler can likewise be applied directly to the target $P(\theta)$ \eqref{eq:P} using exact determinant evaluations.}

Finally, we consider a problem \REV{using real data} and demonstrate convergence of
our method within a tractable amount of time on commercial hardware.

Our full implementation, including code for all our experiments, is
available at \cite{Kielstra_2024}.  \REV{We ran all our computations on a single virtual machine on the Jetstream2 cloud, equipped with AMD Milan 7713 CPUs and an NVIDIA A100 GPU.}

\subsection{\label{subsec:simple}Verification in a simple case}

\REV{We} take $N=10$ points, $N_{\textrm{cheb}}=2$, and $d=1$. \REV{We} also fix
$2\ell^2=1$ and $\sigma_{\theta}^{2}\equiv\sigma^{2}=0.1$. For repeatability
and simplicity, we do not take the scattered points $x_{i}$ to be
random; rather, we let them be equally spaced on the interval $[-1,1)$.
Similarly, we take each observation to be $y_{i}=1$ deterministically.
Since $\theta\in\R^{2}$ in this case, we can numerically compute
$P(\theta)$ on a grid of $100^{2}$ equispaced points on $[-3,3]^{2}$
and perform integrals using the trapezoidal rule. (These particular
numbers were observed to be large enough that the values of our integrations did not meaningfully change when we increased them.) Then, we run each of the following sampling routines,
batching $B=500$ chains in parallel:
\REV{
\begin{enumerate}
\item Random Walk Metropolis (RWM) with a step size of $\Delta t = 0.25$ for $15000$ steps;
\item Hamiltonian Monte Carlo (HMC), with a leapfrog integrator, a step
size of $\Delta t =  0.4$, and $N_{\textrm{int}}=3$, for $5000$ steps; and
\item HMC, with an implicit midpoint (Gauss-Legendre) integrator, a step
size of $\Delta t = 0.15$, and $N_{\textrm{int}}=3$, for $5000$ steps.
\end{enumerate}
For all linear solves as well as the nonlinear solves required by the implicit integrator (cf. Sections~\ref{subsec:varphi}, \ref{subsec:Forces}, and \ref{subsec:implicit}), we use a rank-$5$ randomized Nystr\"om preconditioner, following~\cite{frangella_tropp_udell}. After every two Metropolis-Hastings accept-reject updates, we recompute the preconditioner using a low-rank factorization of $K(\theta)$ for the most recent sample of $\theta$.  We use this preconditioner, unchanged, for the entirety of the next two Metropolis-Hastings updates. Each of the $B$ chains is preconditioned independently of the others.}

\REV{The step sizes were chosen in order to achieve an overall acceptance probability of approximately $0.65$ \cite{Neal_2011}, except in the case of the implicit integrator, for which a smaller step size $\Delta t$ was necessary in order to ensure that the nonlinear equations could be solved.  The preconditioner had little overall effect in the leapfrog case, but sped up convergence for the nonlinear solves in the implicit case by a factor between 2 and 3.  It did not significantly increase the step size $\Delta t$ for which a solution of the nonlinear equations could be attained.}

\begin{figure}
\includegraphics[width=1\textwidth]{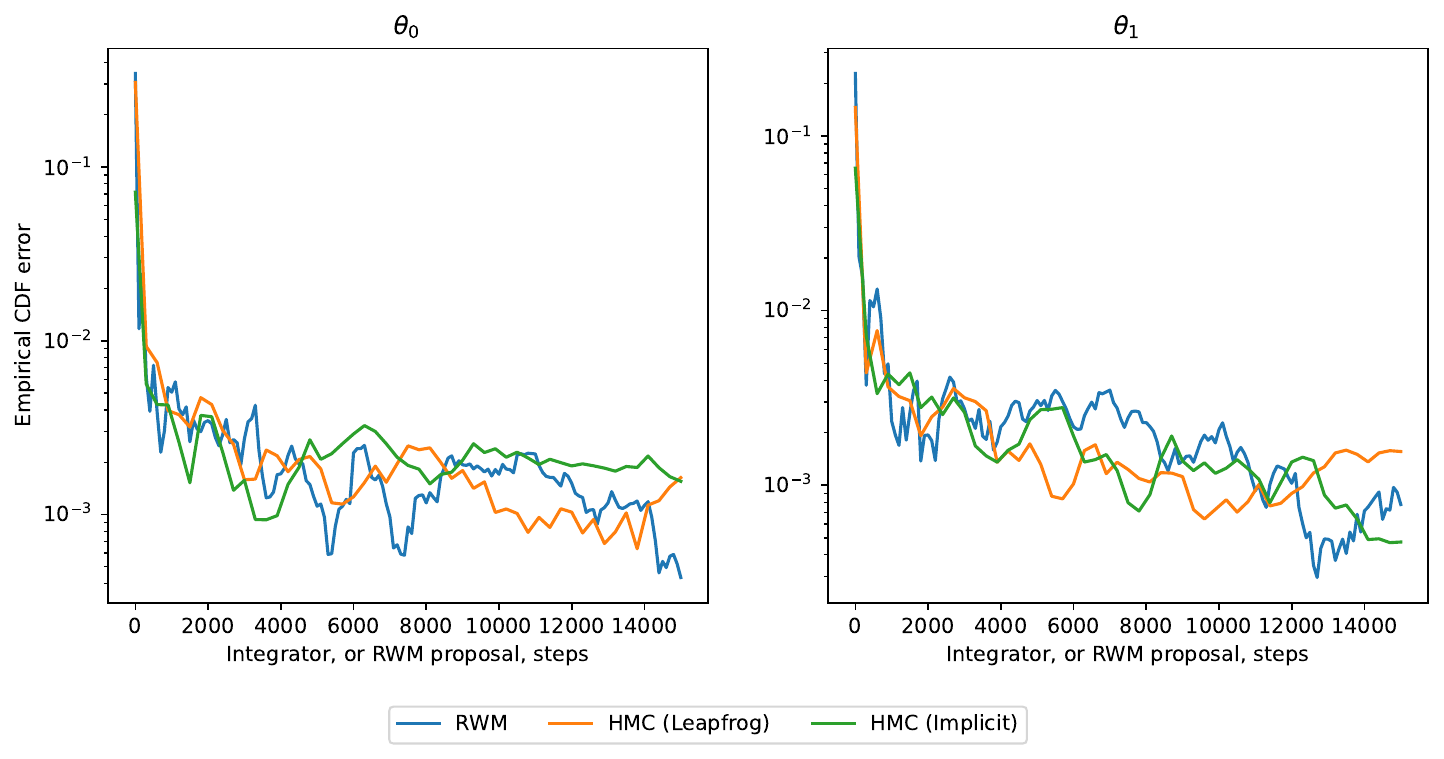}

\caption{Uniform norm error of the empirical marginal CDF for the two components
$\theta_{0}$ and $\theta_{1}$ of the hyperparameter vector $\theta$,
computed with pseudofermion-based sampling routines.\label{fig:theta_marginal} }

\end{figure}

Then, at every step $i$, we computed the empirical
CDF for the marginal distributions of the components $\theta_{0}$
and $\theta_{1}$ of $\theta$ using all $B=500$ batch elements and
samples $\lfloor i/2 \rfloor$ through $i$. This ``moving-window'' procedure \REV{systematically removes burn-in
artifacts}. We compared the results to \REV{the CDF obtained by highly accurate numerical integration} at $100$
\REV{equally-spaced} points in $[-3,3]$, and took the largest difference. The results,
which show good agreement overall, are shown in Figure \ref{fig:theta_marginal}.  \REV{The horizontal axes measure integrator steps, not Metropolis-Hastings proposals, in order to provide for a fairer comparison between HMC, which computes multiple interim hyperparameter values per Metropolis-Hastings proposal, and RWM, which does not.}

% We also include a $95\%$ confidence radius, calculated as per the
% DKW inequality \cite{10.1214/aop/1176990746} as
% \[
% \sqrt{\frac{\log\frac{2}{0.95}}{2\cdot Bi/\tau}},
% \]
% where $i$ is the number of steps taken up to the point where we are
% calculating the bounds and $\tau$ is the empirically-computed \emph{integrated
% autocorrelation time} (IAT) \cite{Goodman_Weare_2010}. This is a measure of the number of samples
% from the chain required to take one unbiased sample of the underlying
% distribution. (See Appendix \ref{app:iat} for more information.)
% We expect $0$ to be within these bounds with probability $0.95$,
% which appears to be compatible with our results. We start the figures
% at $i=50$ rather than $i=0$, as otherwise most of the vertical space
% on the graph would be taken up by the (very short) time taken to reach
% the typical set, and resolution would be lost. 

% Since we ran both the RWM and the HMC routines for the same total
% number of integration or integration-like steps, but the resulting
% RWM chain was three times longer than the HMC chain, we normalize
% by representing each RWM step as being $1/3$ the length of each HMC
% step. Figures \ref{fig:theta_marginal} and \ref{fig:theta_mean}
% both have a fixed ``time'' axis of length $5000$, with a new HMC
% step when time advances by $1$ and a new RWM step when it advances
% by $1/3$.

\begin{figure}

\includegraphics[width=1\textwidth]{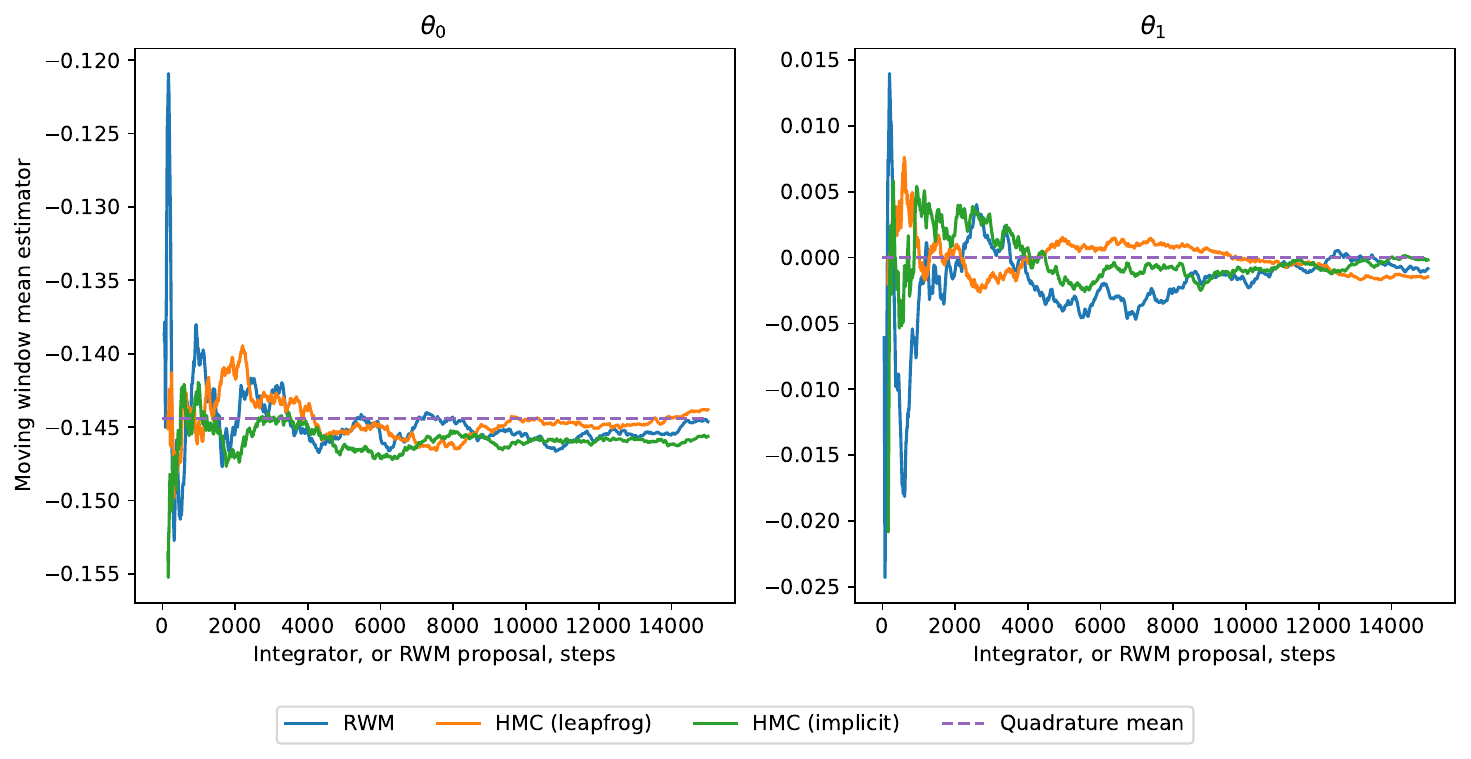}\caption{Error of the empirical means for the two components $\theta_{0}$
and $\theta_{1}$ of the hyperparameter vector $\theta$, computed
with pseudofermion-based sampling routines. \label{fig:theta_mean}}

\end{figure}

To further confirm agreement with the true distribution, we compute
the posterior means for both components of $\theta$ with the trapezoidal
rule. First we average over the batch dimension, and then we use our
samples to estimate the posterior means using a moving-window average.
Specifically, given a chain of samples $X_{1},\dots,X_{n}$, the moving-window
average is a chain $T_{1},\dots,T_{n}$ such that
\[
T_{i}=\frac{X_{\lfloor i/2 \rfloor }+X_{\lfloor i/2 \rfloor +1}+\cdots+X_{i}}{i - \lfloor i/2 \rfloor}.
\]
We plot this chain as computed using RWM, leapfrog HMC, and implicit
HMC and show the results in Figure \ref{fig:theta_mean}. \REV{As in Figure \ref{fig:theta_marginal}, the horizontal axes measure integrator steps, not Metropolis-Hastings proposals.}

We approximate the standard deviation of the estimator $T_{i}$
as
\[
\frac{\sqrt{\textrm{Var}\, X}}{\sqrt{Bi/(2\tau)}},
\]
\REV{where we estimate the variance empirically using chain elements $X_{\lfloor i/2 \rfloor}$ through $X_i$, $B=500$ is the number of parallel chains, and $\tau$ is the empirically-computed \emph{integrated
autocorrelation time} (IAT) \cite{Goodman_Weare_2010}. This is a measure of the number of steps
in the chain required to furnish one independent sample from the target. (See Appendix \ref{app:iat} for more information.)  The results are shown in Table \ref{tab:variances}.}

\begin{table}
    \centering
    \begin{tabular}{l|c c}
    Proposal & $\theta_0$ & $\theta_1$ \\
    \hline
    RWM & $0.00131$ & $0.00174$ \\
    HMC (leapfrog) & $0.00077$ & $0.00065$ \\
    HMC (implicit) & $0.00165$ & $0.00225$ \\
    \end{tabular}   
    \label{tab:variances}
    \caption{Approximate standard deviation of mean estimator for both kernel hyperparameters as computed with pseudofermion-based sampling routines.}
\end{table}

% \REV{In general, we see that the leapfrog HMC converges most quickly, but all three methods appear to have somewhat similar rates of convergence to each other.  In the interests of providing experiments that can be cheaply replicated, we do not run any of the chains to full convergence; in other experiments, we found that the error continued to decrease, albeit more slowly, as the sampling time increased.}

\begin{table}
    \centering
\begin{tabular}{l|c c c}
Proposal & $\theta_0$ & $\theta_1$ \\
\hline
RWM & 1.15e-03 & 1.45e-03 \\
HMC (leapfrog) & 3.94e-04 & 2.02e-04 \\
HMC (implicit) & 1.83e-03 & 2.39e-03 \\
\end{tabular}
    \caption{Time, in seconds, to compute a single unbiased sample, based on an a posteriori calculation of the IAT.}
    \label{tab:unbiased_times}
\end{table}

\REV{We also computed the total wall time per independent sample, defined as $\frac{t_{\mathrm{wall}}}{BL/\tau}$, where $t_{\mathrm{wall}}$ is the total wall time, $B$ the number of parallel chains, $L$ is the length of each chain, and $\tau$ the IAT.  The results are shown in Table \ref{tab:unbiased_times}, and show a clear advantage for leapfrog-based HMC. Note that this is only a 2-dimensional hyperparameter sampling problem, and we expect the advantage of gradient-based methods to only increase in higher dimensions~\cite{Chewi_2024}.}

% We computed the IAT from the full chain, and the variance from its
% full second half, but we only show the uncertainty starting from the
% fiftieth timestep since otherwise it is too large to allow the rest
% of the graph to render readably. Figure \ref{fig:theta_mean}, which
% includes error bars for one standard deviation, shows that the variance
% is much reduced by the use of HMC, and in particular by the use of
% implicit HMC. In general, convergence of RWM is much worse than that
% of HMC.

\subsection{Scaling analysis}

For a more computationally intensive problem, we choose our points
uniformly at random in $[-1,1]^{d}$. We take our observations to be
\[
y_{i}=\prod_{j=1}^{d}\cos(x_{i}^{j})+\tilde{\epsilon}_{i},
\]
where $i=1,\dots,N$, $\tilde{\epsilon}_{i}\sim\mathcal{N}(0,\eta^{2})$ are i.i.d.
and $\eta^{2}=0.01$. \REV{Here, $x_i^j$ is the $j$-th component of the vector $x_i$.} We introduce the tilde notation to emphasize
that the data-generating process may not coincide with the modeling
assumption (\ref{eq:model}) used for inference.

Keeping $N_{\textrm{cheb}}=2$, we run our three samplers \REV{for the determinant-based target~\eqref{eq:P} and the pseudofermion-based target~\eqref{eq:extendeddensity} with a batch size of $B=1$},
for $200$ steps, in dimensions $d=1$, $2$, and $3$. \REV{As before,
we take $N_{\textrm{int}}=3$ for both HMC samplers.  For all samplers, we set the step size to be $\Delta t = 0.01$.}  To slow the growth of the kernel matrix condition number,
we narrow the kernel function as $N$ grows by setting $2\ell^2=\left(\frac{N}{10^{4}}\right)^{-2/d}$
and we take $\sigma_{\theta}^{2}\equiv\sigma^{2}=0.1$. \REV{A narrower kernel function implies less interaction between any two given data points, which in turn produce rows and columns of the kernel matrix that are closer to orthogonal, reducing the kernel matrix condition number overall, cf.~\cite{kielstra2024gaussianprocessregressionloglinear}. The large value of $\sigma_\theta^2$} assumes
more noise than the true value of $\eta$ would suggest, but improves
the conditioning of the kernel matrix, leading to faster solves and
therefore more efficient profiling. The power of $-2/d$ comes from
our desire for the volume of the numerical support of the kernel function
to shrink with $N^{-1/d}$, guaranteeing that the expected number
of points within that volume remains constant.\footnote{In fact, the maximum number of points within such a neighborhood of
one point grows as $O(\log N)$, so the conditioning in this regime
does as well. The effect is small enough to be invisible in our results.} \REV{For the implicit HMC proposal,} we use the simple \REV{rescaling} preconditioner from Section \ref{subsec:implicit}, estimating $c \coloneqq \rho(A)$ with the power method and setting our preconditioner to be $\frac{1}{c}I$.  \REV{We recompute this rescaling factor after every integrator step. Otherwise, no preconditioning is used for any linear solves.} 

\begin{figure}
\includegraphics[width=1\textwidth]{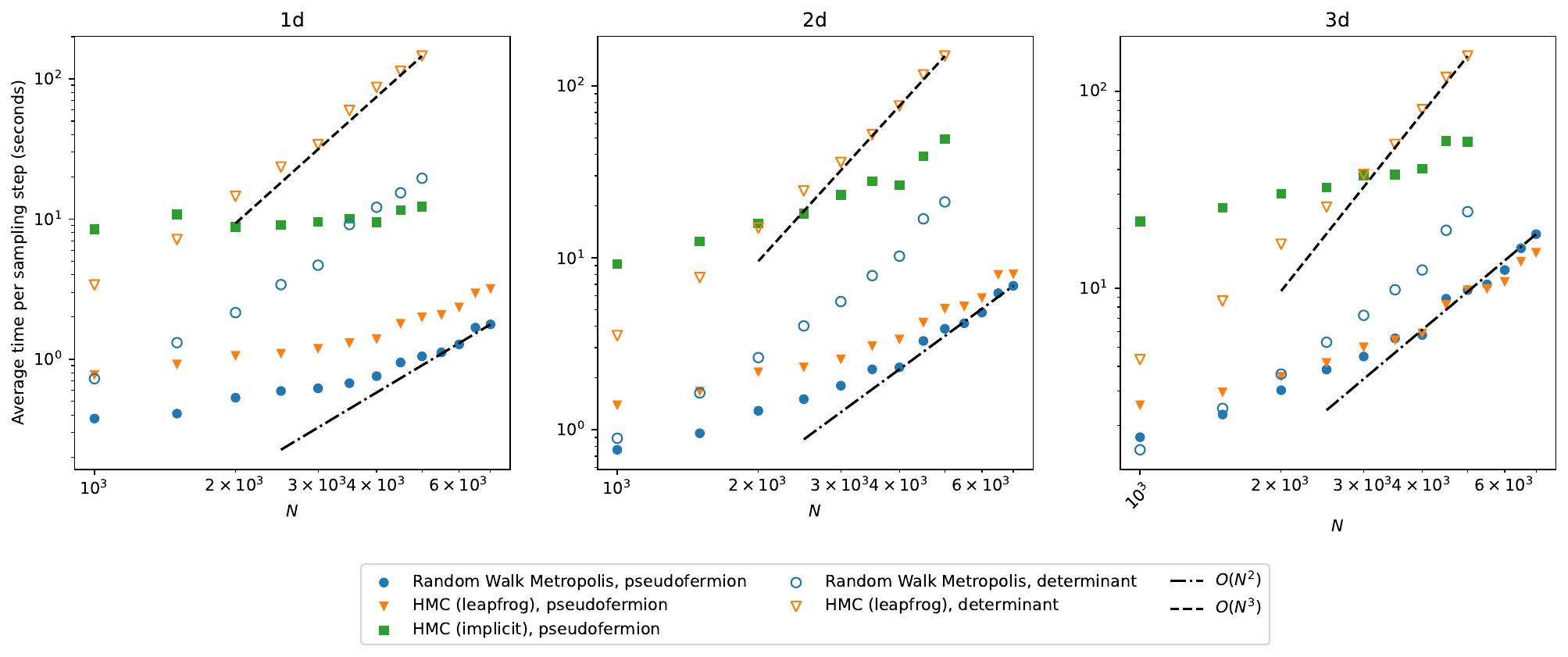}

\caption{Computational time per sampling step of various methods on an A100
GPU. $N$ is the number of data points. \label{fig:scaling} }
\end{figure}

The results are shown in Figure \ref{fig:scaling}. We see
$O(N^{2})$ scaling for our method and $O(N^{3})$ scaling in the
determinant case, as expected. We also see that our method requires
a fairly large value of $N$ before such scaling occurs. This is due
in part to our use of KeOps \cite{JMLR:v22:20-275}, a library which
exploits parallelism in kernel matrix-vector multiplications. Note that the memory footprint
for these operations in KeOps is only $O(N)$, and the entire kernel
matrix need not ever be formed in memory. Determinant-based methods,
meanwhile, do require formation of the whole matrix, leading to increased
memory use and faster accelerator saturation.

In theory, it is possible to use implicit HMC with determinant-based
sampling. In practice, the performance was not competitive enough
to profile. Also in theory, we could have profiled the determinant-based
methods and the implicit pseudofermion methods at the larger values
of $N$ at which we profiled the leapfrog and RWM pseudofermion methods.
In practice, once again, the performance was not competitive enough.

\begin{figure}
\includegraphics[width=1\textwidth]{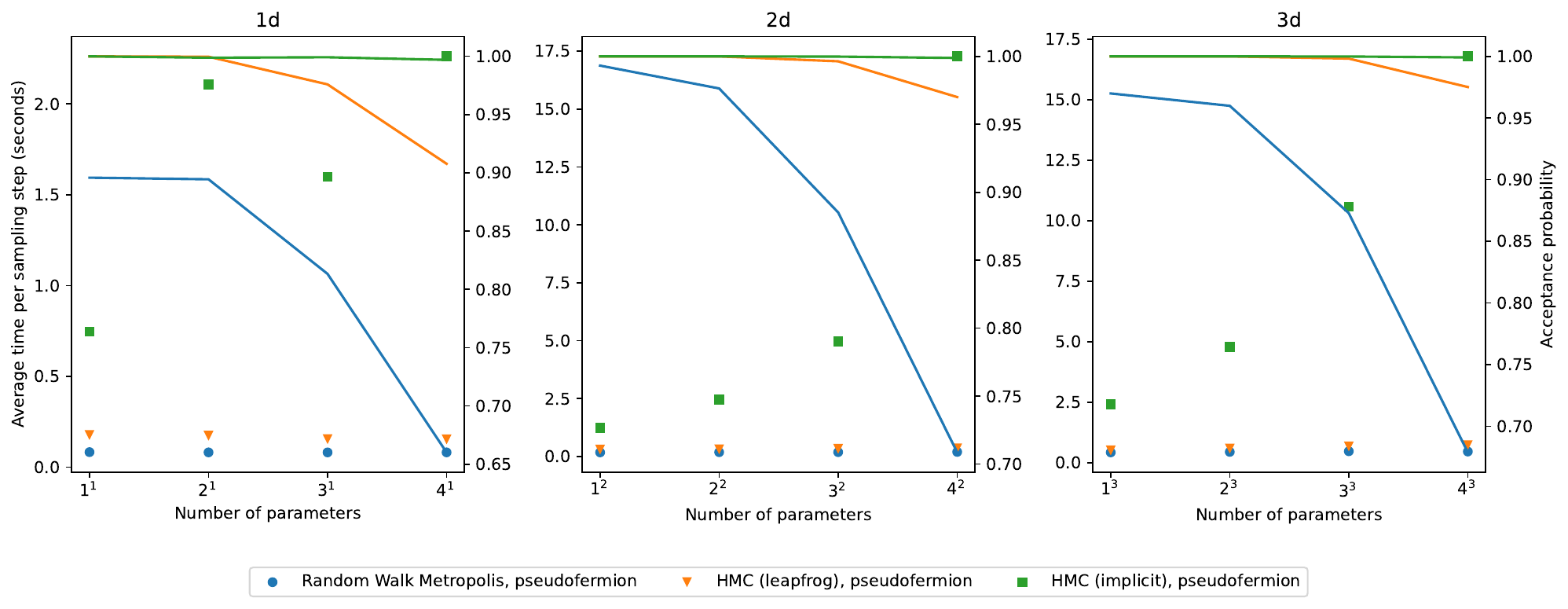}

\caption{Acceptance probability (solid lines) and wall clock (dots) scaling
with number of hyperparameters on A100 GPU.\label{fig:paramscaling}}
\end{figure}

We also fix $N=1500$ and experiment with increasing the Chebyshev
degree $N_{\textrm{cheb}}$ in $d=1,2,3$. The results are shown in
Figure \ref{fig:paramscaling}. Since we take a fixed value of $\Delta t =0.01$,
the Metropolis-Hastings acceptance probability decreases for RWM
and the leapfrog integrator: the increased number of parameters increases
the possibility that a randomly-chosen move will shift away from the
typical set in one or more of them. However, it remains stable for
the implicit integrator. That said, the cost per sampling step for
the implicit integrator increases with the number of hyperparameters
due to the greater number of nonlinear equations to be solved, while
the wall clock scaling of RWM and the leapfrog integrator with respect
to the number of hyperparameters is so mild that it appears empirically
constant in these experiments.

While the implicit integrator does give \REV{large} acceptance probabilities
in this case, this does not translate into the possibility of taking
steps much larger than $\Delta t=0.01$, which can break the numerical convergence
of the accelerated fixed-point iteration required by the numerical
integrator. \REV{The potential computational advantage of the implicit method in even higher-dimensional settings, where the consequences of falling acceptance probability for other methods become more exaggerated, is an interesting topic for further exploration.}

\subsection{A real GPR problem\label{subsec:A-large-GPR}}

Finally, we test our algorithm on real data.  Inspired by \cite{Cressie_2018}, we carry out GPR on satellite measurements of atmospheric CO\textsubscript{2} from \cite{SCF_2017}.  After restricting to a box roughly surrounding the contiguous United States (latitude between $20$ and $50$ and longitude between $-125$ and $-65$), we are left with approximately $110000$ data points, many of which are tightly clustered geographically.  As a preprocessing step, we de-mean the data, linearly shifting the coordinate points so that their mean is at $(0, 0)$ and the CO\textsubscript{2} measurements so that their mean is $0$ as well.  Then we rescale the coordinates into $[-1, 1]^2$ and randomly select $N=3000$ training points.  We reverse the scaling and shifting before plotting the results, but after calculating the posterior GP variances and estimator uncertainties.

We use the same kernel as before, carrying out experiments for $N_{\textrm{cheb}}=1,2,3,4$.  We also infer the values of $\sigma$ and $\ell$.  Specifically, in order to ensure that they are always strictly greater than zero, we introduce hyperparameters $\tilde{\sigma}, \tilde{\ell} \in \mathbb{R}$ and set $$\sigma = \exp(\tilde{\sigma}) + 10^{-3},\quad 2\ell^2 = \exp(\tilde{\ell}) + 10^{-3}.$$  In general, therefore, the number of kernel hyperparameters is given by $N_\textrm{cheb}^2+2$.

\begin{figure}
    \centering
    \includegraphics[width=\linewidth]{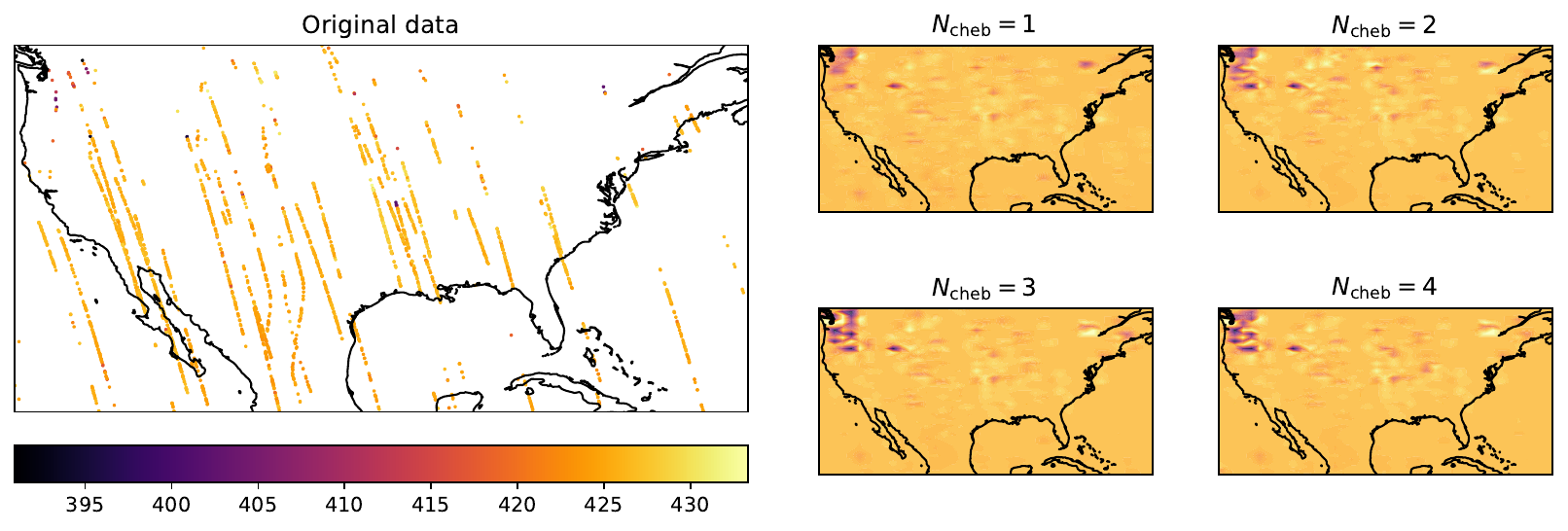}
    \caption{Inferred GPR means for satellite data.}
    \label{fig:satellite}
\end{figure}

For each value of $N_{\mathrm{cheb}}$, we run a batch of $B=10$ chains in parallel, using HMC with the leapfrog integrator and $N_\textrm{int}=3$.  We use the same basic preconditioning strategy as in Section \ref{subsec:simple}, recalculating a rank-500 Nystr\"om preconditioner after every five Metropolis-Hastings updates.  We also use this preconditioner to accelerate the sampling of $\varphi$, as detailed in Section \ref{subsec:varphi}.  We first run a burn-in procedure consisting of $300$ updates with a step size of $\Delta t = 0.01$, followed by $12500$ updates with a step size of $\Delta t = 0.03$.  (The increased step size becomes more efficient as we approach the typical set.)  Then, we discard the burn-in samples and run the chain for $12500$ more updates, again with a step size of $\Delta t= 0.03$.  The averaged Metropolis-Hastings acceptance probabilities are shown in Table \ref{tab:bigprobs} and fall between 75\% and 85\%.  For each hyperparameter sample, we calculate the inferred GP mean function (with kernel determined by this hyperparameter sample) at every point on a $50 \times 50$ grid, following~\cite{Rasmussen_Williams_2006}.  Then we empirically average over the hyperparameter uncertainty using our samples and plot the resulting function in Figure \ref{fig:satellite}. 

\begin{table}[]
    \centering
    \begin{tabular}{c|c}
$N_\textrm{cheb}$ & Mean acceptance probability \\
\hline
1 & 81.09\% \\
2 & 79.96\% \\
3 & 78.31\% \\
4 & 77.17\% \\
\end{tabular}
    \caption{Average Metropolis-Hastings acceptance probabilities for this experiment.}
    \label{tab:bigprobs}
\end{table}

\begin{table}[]
    \centering
    \begin{tabular}{c|c c c c c c}
$N_\textrm{cheb}$ & Range of values & \makecell{Mean \\ IAT} & \makecell{Max \\ IAT} & \makecell{Mean 95\% \\ CI radius} & \makecell{Max 95\%  \\ CI radius} \\
\hline
(Original data) & $[390.91, 433.34]$ \\
1 & $[402.76, 434.30]$ & 29.61 & 330.08 & 4.83e-05 & 2.83e-03 \\
2 & $[381.94, 454.98]$ & 66.37 & 304.40 & 4.58e-04 & 3.10e-02 \\
3 & $[374.50, 545.31]$ & 139.92 & 394.95 & 2.31e-03 & 1.16e+00 \\
4 & $[362.93, 463.73]$ & 383.25 & 1101.29 & 2.88e-03 & 4.06e-01 \\    \end{tabular}
    \caption{Analysis of GPR mean inference, including ranges of final estimated values and mean and maximum values for IATs and for 95\% confidence interval radii for those final values.}
    \label{tab:bigresults}
\end{table}

In order to make the plots of the GPR results easier to read, they were clipped to lie within the range of values of the original data.  The pre-clipping ranges are shown in Table \ref{tab:bigresults}, along with the results of a convergence analysis.  We consider the chain of inferred values for every point on the grid (averaged over the batch dimension) and compute its IAT.  Then, again considering every grid point separately, we compute an approximate 95\% confidence interval for our estimation of the GP mean. Similarly to Section \ref{subsec:simple}, we estimate the radius of this interval as $2\frac{\sqrt{\textrm{Var} \,  X}}{\sqrt{12500B/\tau}}$, where $X$ is the inferred GP mean at a given grid point and $\tau$ is the IAT.  We estimate the variance and IAT empirically using all 12500 hyperparameter samples.  The factor of $2$ doubles the standard deviation to give the radius of a $95\%$ confidence interval.  Note that this is not an estimate of the posterior GP variance; it is a measure of our numerical error in computing the expectation of the GP mean function with respect to the hyperparameter uncertainty. 
% Due to the de-meaning procedure mentioned earlier, the estimated means cluster around $0$, so, while the confidence-interval radii are still valid after shifting back to the original mean, they should not be used to compute relative errors based on that specific mean value. 
We report the mean and maximum values of both the IAT and the confidence interval radius across all grid points.

We also compute an estimate of the total uncertainty in the inferred function, combining the contributions due to the GP uncertainty (for fixed kernel hyperparameters) and the uncertainty of the hyperparameters themselves.  Specifically, let $f : \R^2 \ra \R$ denote the unknown function to be inferred, which we assume to be GP-distributed according to the kernel $\mathcal{K}_{\theta}$, conditioned on kernel hyperparameters $\theta$.  For any $x \in [-1, 1]^2$, the law of total variance then allows us to compute $$\textrm{Var}( f(x) ) = \mathbb{E}[\textrm{Var}( f(x) \mid \theta)] + \textrm{Var}(\mathbb{E}[ f(x) \mid \theta]).$$  The second term can be estimated by computing the inferred GP mean for each sample of the hyperparameters (as in the preceding discussion) and then estimating the empirical variance of this quantity over our hyperparameter samples. 

To calculate the first term, let $K_\theta (x) \in \mathbb{R}^N$ be the vector with $j$-th element given by $\mathcal{K}_\theta(x, x_j)$.  Then, given our observations and conditioning on $\theta$, the variance at $x$ is given by~\cite{Rasmussen_Williams_2006} 
\begin{equation}\label{eq:condvar}
\textrm{Var}( f(x) \mid \theta) =  \mathcal{K}_{\theta} (x,x) - K_\theta (x)^\top A(\theta)^{-1}K_\theta (x).
\end{equation}
This term can then be estimated by computing the conditional variance for each sampled vector $\theta$ of hyperparameter values using~\eqref{eq:condvar}, then empirically averaging over samples.

\begin{figure}
    \centering
    \includegraphics[width=\linewidth]{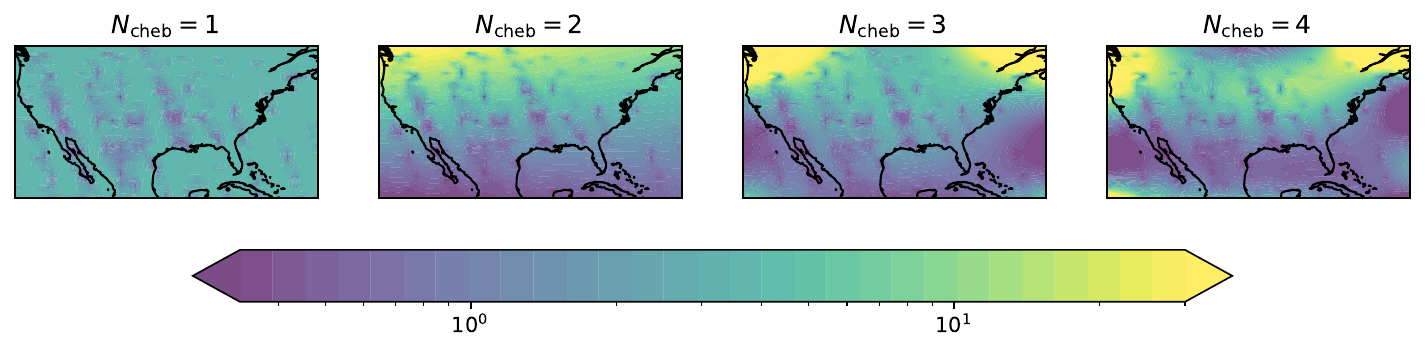}
    \caption{Inferred GPR standard deviations for satellite data.}
    \label{fig:uncertainties}
\end{figure}

We calculate the standard deviation $\sqrt{\mathrm{Var}(f(x))}$ at every point $x$ in our grid.  In order to avoid redundant computation due to autocorrelation in the sampler, when averaging over hyperparameter samples, we use only every $125$-th sample from our chains of length $12500$, yielding a total of $100$ samples for each of our $B=10$ parallel chains, or $1000$ samples in total. We plot the results in Figure \ref{fig:uncertainties}.  In the case $N_\textrm{cheb}=1$, the uncertainty is bounded above by $\approx 13$.  However, as $N_\textrm{cheb}$ increases, the posterior variance grows in places with less data.  For this reason, we clip all posterior standard deviations to the range $[1/3, 30]$ before plotting; the brightest areas on the graph should be read as having insufficient data for meaningful inference.

For each value of $N_\textrm{cheb}$, the experiment took less than four hours. We imagine that such experiments could be accelerated even further by incorporating more advanced methods to apply, or invert, the kernel matrix. Indeed, our dense $O(N^2)$-scaling matrix-vector multiplications dominate the computational cost, but for certain kernels this operation can be approximated with just $O(N \log N)$ operations \cite{greengard_equispaced_2023, kielstra2024gaussianprocessregressionloglinear}. That said, our existing implementation is extremely flexible, exploits GPU parallelism naturally, and has only a linear-scaling memory footprint. We also note that the number of hyperparameters does not have a material effect on the cost per update in these experiments. Therefore we conclude that our approach offers a scalable and flexible tool for fully Bayesian Gaussian process regression involving many unknown kernel hyperparameters.

\bibliographystyle{plain}
\bibliography{gpr_pseudofermion}

\appendix

\section{Pole expansion for the inverse matrix square root \label{app:hht}}

Following \cite{doi:10.1137/070700607}, we note that, for an analytic
matrix function $f(A)$, a generalized form of the Cauchy integral
formula allows us to write
\[
f(A)=\frac{A}{2\pi i}\int_{\Gamma}z^{-1}f(z)(zI-A)^{-1}\,dz,
\]
or, fixing some vector $v$, 
\[
f(A)v=\frac{A}{2\pi i}\int_{\Gamma}z^{-1}f(z)(zI-A)^{-1}v\,dz,
\]
where $\Gamma$ is a closed contour containing all the eigenvalues
of $A$. The effect, once this
integral is discretized, is to reduce the problem of computing $f(A)v$
to that of computing a finite number of linear solves of the form
$zI-A$.

For concreteness, we now describe the specialized details of a contour
integral formulation for the matrix function $f(A)=A^{1/2}$, following
\cite{doi:10.1137/070700607}.

Let $M$ be an upper bound on the largest eigenvalue of $A$, which
we estimate with a power method (using ten iterations in the numerical
experiments) and $m$ a lower bound on the smallest eigenvalue. Given
the form of $A$ as specified in (\ref{eq:kerdef}), we always take
$m=\min_{i}\sigma_{\theta}(x_{i})$. Define $k=\sqrt{m/M}$, and let
$w(t)=\sqrt{m}\cdot\textrm{sn}(t|k^{2})$ be a helper function derived
from a particular Jacobi elliptic function. Similarly, write $\textrm{cn}(t)=\textrm{cn}(t|k^{2})$
and $\textrm{dn}(t)=\textrm{dn}(t|k^{2})$. Let 
\[
K=K(k^{2}),\quad K'=K(1-k^{2})
\]
be specific values of the complete elliptic integral of the first
kind. Then, in the specific case $f(A)=A^{1/2}$, the $N_{p}$-point
discretization
\begin{equation}
A^{1/2}\approx f_{N_{p}}(A)\coloneqq\frac{2K'm^{1/2}A}{\pi N_{p}}\sum_{j=1}^{N}\textrm{cn}(t_{j})\textrm{dn}(t_{j})(A-w(t_{j})^{2}I)^{-1},\label{eq:msqdisc}
\end{equation}
where 
\[
t_{j}=i\left(j-\frac{1}{2}\right)K'/N_{p},
\]
converges at the rate
\[
\Vert A^{1/2}-f_{N_{p}}(A)\Vert=O\left(e^{\epsilon-2\pi KN_{p}/K'}\right)
\]
for any $\epsilon>0$ \cite{doi:10.1137/070700607}. Furthermore,
the constant $K/K'$ is only $O(1/\log(M/m))$, while a naive approach
using equispaced integration points would lead to convergence with
an $O(1/(M/m))$ constant in the exponential instead. The result is
a highly accurate integration even with a very small value of $N_{p}$.
In our numerical experiments, we use only $N_{p}=15$.

Since the indices $t_{j}$ are purely imaginary, the shifts $-w(t_{j})^{2}$
are real and positive and the coefficients for the solves are real.
Therefore, as detailed in \cite{doi:10.1137/070700607}, both the
weights and the shifts can be computed with nothing more than SciPy's
\texttt{ellipk} and \texttt{ellipj} functions \cite{2020SciPy-NMeth}
and the solves themselves can be calculated with a batched conjugate
gradient routine.

To compute $A^{-1/2}=A^{-1}A^{1/2}$, we simply remove the first $A$
in (\ref{eq:msqdisc}).

\section{Anderson acceleration \label{app:aa}}

The general fixed-point problem
\[
x=g(x)
\]
can also be written
\[
x-g(x)=0.
\]
If we define the residual 
\[
f(x)=x-g(x),
\]
the problem becomes one of minimizing $\Vert f(x)\Vert^{2}$.

In pure fixed-point iteration, this is done by taking a sequence of
points $x_{n}=g(x_{n-1})$ such that $f(x_{n})\ra0$ as $n\ra\infty$.
By contrast, Anderson acceleration \REV{\cite{Walker_Ni_2011}} uses the history of the previous
$k$ iterates in the sequence to determine the next iterate at each
stage, where $k$ is a user-defined parameter.  \REV{In our implementation, we take $k=10$.}

Specifically suppose we are given iterates $x_{1},\ldots,x_{n-1}$
and seek to define the next iterate $x_{n}$. We choose coefficients
$\alpha_{1},\ldots,\alpha_{k}$ to minimize 

\[
\left\Vert \sum_{i=1}^{k}\alpha_{i}f(x_{n-i})\right\Vert ^{2},\quad\text{subject to }\ \sum_{i=1}^{k}\alpha_{i}=1.
\]
 This can be solved tractably as a constrained linear least squares
problem. Note that we do not require any additional evaluations of
$f$ (or equivalently, $g$) beyond those already performed once on
each iterate. Then the next iterate is defined in terms of the solution
$(\alpha_{i})$ as 
\[
x_{n}=\sum_{i=1}^{k}\alpha_{i}g(x_{n-i}).
\]

It is demonstrated in \cite{Walker_Ni_2011} that, if we set $g(x)=(I-A)x+b$
and carry out Anderson-accelerated fixed-point iteration with full
history (i.e., $k=n$ at the $n$-th iteration) to define a sequence
$\{x_{n}\}$ and simultaneously apply non-truncated GMRES to $Ay=b$
to get a sequence $\{y_{n}\}$, then, as long as $A$ is such that
the GMRES residual norm decreases at every step, then $x_{n+1}=g(y_{n})$.
This motivates our hybrid use of Anderson acceleration to solve the
linear and nonlinear components of our problem simultaneously: if
the nonlinear component converges quickly for a given approximation
to $A^{-1}b$, then the convergence overall should be about as good
as that of the GMRES inversion of $A$. However, note that in our
implementation, we make no particular distinction between different
function arguments for acceleration purposes. All arguments are combined
into a single vector, to which Anderson acceleration is then applied.

\section{The integrated autocorrelation time (IAT)\label{app:iat}}

The accuracy of an unbiased Markov chain estimator is determined by
the behavior of its variance as the chain grows long. In particular,
if we want to estimate
\[
A=\mathbb{E}_{\pi}[f(X)],
\]
where $X$ is a random variable with density $\pi$, then, given a
path $(X_{1},\dots,X_{N})$ drawn from a Markov chain with invariant
distribution $\pi$, our estimator for $A$ is 
\[
\hat{A}_{N}=\frac{1}{N}\sum_{t=1}^{N}f(X_{t}).
\]
 We assume that $X_{1}$ is drawn from the invariant distribution
(which in practice is achieved by removing iterates before some burn-in
time) so that the process is stationary.

\REV{To empirically determine the variance, and hence the quality, of our estimator $\hat{A}_{N}$, we follow \cite{Goodman_Weare_2010}.} We first define
\[
C(j)=\textrm{cov}[f(X_{0}),f(X_{j})],\quad j=0,1,2,\ldots.
\]
 and extend by symmetry $C(j)=C(-j)$ to be defined for all integers
$\mathbb{Z}$.

Then by stationarity, 
\[
\textrm{cov}[f(X_{s}),f(X_{t})]=C(t-s)
\]
 for all $t,s=0,1,2,\ldots$

Then we may compute 
\begin{align*}
\textrm{var}[\hat{A}_{N}] & =\frac{1}{N^{2}}\sum_{t,s=1}^{N}C(t-s)\\
 & =\frac{1}{N}\sum_{j=-N+1}^{N-1}\left(1-\frac{|j|}{N}\right)C(j)\\
 & \approx\frac{1}{N}\sum_{t=-\infty}^{\infty}C(t)\\
 & =\frac{1}{N}\left(\sum_{t=-\infty}^{\infty}\frac{C(t)}{C(0)}\right)C(0),
\end{align*}
 where the approximate equality holds for large $N$ (assuming, for
example, exponential decay of $C(t)$ as $\vert t\vert\ra\infty$).

Defining the integrated autocorrelation time (IAT) as 
\[
\tau=\sum_{t=-\infty}^{\infty}\frac{C(t)}{C(0)},
\]
then, for $N$ sufficiently large, 
\[
\textrm{var}[\hat{A}_{N}]\approx\frac{\textrm{var}_{\pi}[f(X)]}{N/\tau},
\]
giving us the intuitive understanding that $N/\tau$ is the number
of ``effectively independent'' samples furnished by a chain of length
$N$.

We use the \texttt{emcee} Python package to compute $\tau$ in practice
\cite{Foreman-Mackey_Hogg_Lang_Goodman_2013}. It uses the estimator
\[
\hat{C}(j)=\frac{1}{N-j}\sum_{i=1}^{N-j}(f(X_{i})-\hat{A}_{N})(f(X_{i+j})-\hat{A}_{N}),\quad j=0,\ldots,N-1.
\]
Then, rather than estimate the IAT as 
\[
\sum_{j=-N}^{N}\frac{\hat{C}(\vert j\vert)}{\hat{C}(0)}=1+2\sum_{j=1}^{N}\frac{\hat{C}(j)}{\hat{C}(0)},
\]
\texttt{emcee} in fact estimates the IAT as
\[
\hat{\tau}=1+2\sum_{j=1}^{M}\frac{\hat{C}(j)}{\hat{C}(0)}
\]
 for some $M\ll N$. The reason for this is that, due to the small
amount of data available for the estimator to work with, $\hat{C}(j)$
is very noisy when $j$ approaches $N$. It is therefore more effective
to take $M$ as small as possible while still retaining enough terms
in the sum to capture the true IAT. We can be confident that we have
done this if 
\[
M>c\hat{\tau}
\]
 for some fairly large constant $c$. The \texttt{emcee} package takes
$c=5$ by default. Finally, in order to ensure that the resulting
$M$ is sufficiently small compared to $N$, \texttt{emcee} by default
requires that $N\geq50\hat{\tau}$. \REV{This condition is} easily satisfied
in our \REV{first experiment}.  \REV{In our final experiment, we are only able to assert that $N \geq 25\hat{\tau}$ in most cases, and for the most complicated sampling we achieve little more than $N \geq 10\hat{\tau}$, but these values are still large enough that, for the purposes of IAT estimation, we did not see much need to run a longer chain.}
\end{document}